%% file: root.tex
\documentclass[letterpaper, 10 pt, conference]{ieeeconf}  

\IEEEoverridecommandlockouts                              

\usepackage{multirow}
\usepackage{cite}
\input{preamble}

\title{\LARGE \bf
  Constrained Policy Optimization for Stochastic Optimal Control under Nonstationary Uncertainties*
}

\author{
  Sungho Shin,$^{1}$
  François Pacaud,$^{1}$
  Emil Contantinescu,$^{1}$
  Mihai Anitescu$^{1,2}$
  \thanks{*This material is based upon work supported by the U.S. Department of Energy, Office of Science, Office of Advanced Scientific Computing Research (ASCR) under Contract DE-AC02-06CH11347.}
  \thanks{$^{1}$Mathematics and Computer Science Division, Argonne National Laboratory, Lemont, IL 60439}%
  \thanks{$^{2}$Department of Statistics, University of Chicago, Chicago, IL 60637}
}

\begin{document}

\maketitle

\begin{abstract}
  This article presents a constrained policy optimization approach for the optimal control of systems under nonstationary uncertainties.
  We introduce an assumption that we call
  \emph{Markov embeddability} 
  that allows us to
  cast the stochastic optimal control problem as a policy optimization problem over
  the augmented state space.
  Then, the infinite-dimensional policy optimization problem is approximated as a finite-dimensional nonlinear program by applying function approximation, deterministic sampling, and temporal truncation.
  The approximated problem is solved by using automatic differentiation and condensed-space interior-point methods.
  We formulate several conceptual and practical open questions regarding the asymptotic exactness of the approximation and the solution strategies for the approximated problem.
  As a proof of concept, we provide a numerical example demonstrating the performance of the control policy obtained by the proposed method.
\end{abstract}

\section{Introduction}

Sequential decision problems under nonstationary uncertainties have been recognized as  computationally challenging \cite{bertsekas2012dynamic}.
The stochastic problem admits a deterministic equivalent for linear-quadratic settings with additive uncertainties (so-called certainty equivalence \cite{theil1957note,simon1956dynamic,jacobson1973optimal}), but this result does not generalize to systems with nonlinearity, constraints, or multiplicative uncertainties.
Numerical solution methods for  stochastic sequential decision problems have been studied extensively in the context of multistage stochastic programming \cite{birge1985decomposition}.
Several decomposition algorithms, such as nested Benders decomposition \cite{ho1974nested}, progressive hedging \cite{mulvey1991applying,rockafellar1991scenarios}, and stochastic dual dynamic programming \cite{pereira1991multi} have been proposed.
These methods have limitations, however, in that they rely heavily on convexity and require expensive iterative computations.

In this paper we address the challenges in stochastic control under nonstationary uncertainty based on the {\it Markov embeddability} assumption.
This condition assumes that the uncertainty that drives the stochastic system is a partial observation of a higher-dimensional Markov process.
This assumption can be justified because, in many applications, the uncertain exogenous factors that drive the system (e.g, energy price and demand) can be forecast by the simulation on higher-dimensional spaces with augmented data (e.g., spatiotemporal climate models), which in principle satisfy the Markov property.
This assumption provides a desirable structure that allows us to apply the state augmentation technique \cite{bertsekas2012dynamic} and to embed the exogenous nonstationary process on a single Markov decision process.
In turn, one can reformulate the original stochastic control problem as a time-independent policy optimization (PO) problem, which can be subsequently approximated by a finite-dimensional nonlinear program using several approximation strategies (function approximation, sample average approximation, and temporal truncation).
The resulting problem is computationally challenging, but a scalable solution is possible by applying structure-exploiting methods---batched automatic differentiation  and condensed-space interior-point methods (IPMs).
These methods are particularly well suited for running on  modern graphics processing unit/single instruction, multiple data (GPU/SIMD) architecture.

Policy optimization approaches are widely used, classically in the context of dynamic programming (DP), and more recently in the context of reinforcement learning (RL) \cite{bertsekas2012dynamic}.
Since the method proposed in this article can be used in conjunction with the system identification methods, our method in principle can be used as a model-based policy optimization method, which is similar in spirit to the model-based RL approaches \cite{janner2019trust}.
However, our approach is different from the algorithms studied in the RL context in the following ways:
(i) we apply  deterministic sampling rather than  stochastic algorithms \cite{bottou2018optimization}; 
(ii) we do not replace the constraints using a penalty, which is standard in RL, but instead apply the classical constrained nonlinear optimization algorithm \cite{nocedal1999numerical}; and
(iii) our method does not include online learning, which is an important aspect of many RL methods.
In the classical control context, several works have studied the policy optimization approach.
The work by Fazel and coworkers established the global convergence property of  zeroth-order and natural gradient algorithms for policy optimization under linear-quadratic settings \cite{fazel2018global}.
Our work is close in spirit to the recent study on the direct policy optimization method applied to trajectory optimization problems \cite{howell2021direct}.
While our approach is similar to that  study in terms of problem formulations and approximation strategies, we undertake a deeper investigation into how nonstationary uncertainties are handled and  what  the practical considerations are in implementing scalable algorithms for policy optimization.

The contributions of this paper are threefold. First, we present  constrained policy optimization based on state augmentation and finite-dimensional approximation strategies.
Second, we highlight several conceptual open questions and computational challenges.
Third, as a proof of concept, we numerically demonstrate the effectiveness of  policy optimization relative to model predictive control (MPC) and linear quadratic regulator (LQR).

\section{Settings}
\subsection{Model}\label{sec:setting-1}
\FP{Let $(\Omega, \cF, \bbP)$ be a probability space.}
We consider a discrete-time {\it stochastic process} $\{\xi_t\}_{t=0}^\infty$, where $t$ is the time index and $\xi_t$ is a random variable taking a value in some measurable set $\Xi\subseteq \mathbb{R}^{n_\xi}$.
\FP{We assume that the distribution of the full uncertainty is known:
at time $t$, the knowledge of the conditional distribution for the future uncertainties $\bxi_{t+1:\infty}$ given $\bxi_{0:t}$ can be exploited.}
The stochastic dynamical system under study is \FP{given by the measurable function
$f: \bbR^{n_x} \times \bbR^{n_u} \times \bbR^{n_\xi} \to \bbR^{n_x}$, such that}
\begin{equation}\label{eqn:dyn}
x_t = f(x_{t-1},u_{t-1};\xi_{t}),
\quad \forall t = 0, \cdots, T-1,
\end{equation}
where $x_t\in\mathbb{R}^{n_x}$ and $u_t\in\mathbb{R}^{n_u}$ denote the state and control variables at time $t$.
We seek to minimize the discounted summation of the {\it performance index} $\ell(x_t,u_t;\xi_t)$ over an infinite horizon while satisfying
\FP{almost surely the} algebraic constraints
\begin{equation}\label{eqn:pi}
  g(x_{t},u_{t},\xi_{t}) \geq 0 ,
\quad \forall t = 0, \cdots, T-1,
\end{equation}
\FP{with $g: \bbR^{n_x} \times \bbR^{n_u} \times \bbR^{n_\xi} \to \bbR^{m}$ a measurable mapping.}
Chance constraints \cite{paulson2020stochastic} are often considered to mitigate the conservativeness of the almost-sure constraints.
Here we assume that the uncertainty distribution does not have a long tail, and thus the problem with almost-sure constraints provides a reasonable control policy.

\FP{According to the dynamics~\eqref{eqn:dyn},}
the following event order is assumed:
\begin{equation*}
  \xi_0, x_0, u_0,\xi_1, x_1, u_1,\cdots,u_{T-1},\xi_T, x_T.
\end{equation*}
In each stage, the control decision $u_t$ is a recourse decision made after partially observing the past uncertainties $\bxi_{0:t}$~\cite{birge2011introduction}:
\FP{
for all $t$, $u_t$ can be expressed by a measurable \emph{mapping} $u_t(\cdot)$ such that
$u_t = u_t(\bxi_{0:t})$.}

An important feature of our setting is that the system under study is driven by a nonstationary stochastic \FP{process} $\{\xi_t\}_{t=0}^\infty$.
\FP{We assume the stochastic process is \emph{non-Markovian},}
in contrast to the classical stochastic control literature where it is typically assumed that $\xi_t$ are mean zero, independent, and identically distributed (iid).
This assumption is proper if the controller's goal is to track a setpoint in the face of mean zero iid disturbances.
However, in many emerging applications such as energy management \cite{chen2022reinforcement}, there is no clear notion of setpoints, and the control goal is rather directly optimizing the economic performance by responding to the time-varying uncertainty signals (e.g., renewable generation and demand).
This is because the system is driven by nonstationary exogenous factors involving time scales comparable to the one of the dynamics \cite{dowling2018economic,tsay2019optimal}.
Remarkably, weather conditions, which have strong correlation with the renewable generation and demand, do not satisfy the Markov property by themselves, since accurately modeling the climate dynamics typically requires the use of long memory.
Our setting is suitable for such applications since we do not make restrictive assumptions about the nature of uncertainty.

\subsection{Problem Formulation}\label{sec:setting-2}
The settings introduced in Section \ref{sec:setting-1} give rise to the following stochastic optimal control problem:
\begin{subequations}\label{eqn:orig}
  \begin{align}
    \min_{\substack{ \{u_t(\cdot)\}_{t=0}^\infty}}\;\label{eqn:orig-obj}
    & \mathbb{E}_{\bxi}\left[\sum_{t=0}^\infty \gamma^t \ell(z_t(\bxi_{0:t}); \xi_t)\right]\\
    \st\;
    &x_{0}(\xi_0)\sim \mathcal{X}\\
    &x_{t}=f(z_{t-1}(\bxi_{0:t-1}); \xi_t),\; a.s.,\;t\in\mathbb{I}_{\geq 1}\label{eqn:orig-dyn}\\
    &g(z_{t}(\bxi_{0:t}); \xi_t) \geq 0,\; a.s.,\;t\in\mathbb{I}_{\geq 0}\label{eqn:orig-con},
  \end{align}
\end{subequations}
where $z_t(\cdot):=(x_t(\cdot),u_t(\cdot))$, $\mathcal{X}$ is the distribution of the initial state, and $\gamma\in (0,1)$ is the discount rate.
We do not express $\{x_t(\cdot)\}_{t=0}^\infty$ as decision variables because they can always be eliminated and expressed in terms of $\{u_t(\cdot)\}_{t=0}^\infty$.
Problem \eqref{eqn:orig} has appeared in various contexts, with different variations.
For example, in the DP and RL literature, the dynamics \eqref{eqn:orig-dyn} is expressed as transition probabilities over (typically) the finite state-action space \cite{bertsekas2012dynamic}.
In the multistage stochastic programming literature, the problem is expressed as the minimization of nested expected values \cite{birge2011introduction}.
The problem \eqref{eqn:orig} is general in that it covers both cases.

Solving \eqref{eqn:orig} is challenging because of its infinite-dimensional nature,
following from (i) the continuous domain of the decision functions $u_t(\cdot)$, (ii) the expectation taken over continuous random variables, and (iii) the infinite time horizon.
Without a simplifying assumption, solving \eqref{eqn:orig} is intractable in general.
The solution approaches for \eqref{eqn:orig} have been of particular interest in the stochastic programming literature.
Sample average approximation  \cite{shapiro1990differential} is the standard technique in dealing with stochastic optimization problems over continuous distributions.
Under suitable assumptions, the solution of a finite-dimensional approximated problem asymptotically converges to that of the original problem \cite{shapiro2009time}.
However, these approaches result in exponentially large scenario trees that need to be explicitly considered within a single optimization problem and thus are not scalable.

\section{Policy Reformulation}
To overcome the computational challenges faced in \eqref{eqn:orig}, we introduce the {\it Markov embeddability} assumption.

\begin{assumption}\label{ass:markovian}
  There exists a Markov \FP{stochastic} process $\{\zeta_{t}\}_{t=0}^\infty$ with $\zeta_t$ being a random variable taking a value in some measurable set $Z\subseteq \mathbb{R}^{n_\zeta}$ and a \FP{measurable} mapping $\psi: Z \to \Xi$ such that $\xi_t = \psi(\zeta_{t})$.
\end{assumption}
In other words, we assume that $\{\xi_t\}_{t=0}^\infty$ is a partial observation of potentially a higher dimensional Markov process.

For the simplest case, $\zeta_t$ may be constructed based on the time-delay embedding \cite{takens1981detecting}. Furthermore, several more sophisticated models can be used to
construct the Markov process, such as the finite-memory
Koopman/Mori--Zwanzig (K-MZ)
formalism~\cite{lin2021data,Chorin_2015a,uy2021operator,williams2015data},
multifidelity model reduction
strategies~\cite{renganathan2022camera,chandy2010t}, and
physics-informed strategies~\cite{Constantinescu_2013b}.
For example, if the probability law of $\bzeta_{0:t}$ is 
known, the K-MZ formalisms allow us to obtain a minimal
representation that trades computational efficiency for accuracy. If,
however, $\bzeta_{0:t}$ is only sparsely observed, then K-MZ provides
us with a basis for inferring the dynamics of the Markov embedding. 

This assumption allows for the construction of an augmented Markov decision process from the stochastic dynamical systems driven by nonstationary uncertainties.
\FP{This is in direct contrast to the classical approaches~\cite{bertsekas1996stochastic,puterman2014markov},
  which render the problem Markovian by introducing a history process $h_t = (x_0, u_0, x_1,\cdots, u_{t-1}, x_t)$, collecting all the past information, and looking
  at control policies $u_t(h_t)$ depending on the current history $h_t$.
  Despite its theoretical interest, this approach is intractable since the dimension of the history variable $h_t$ increases linearly with the time $t$.
  In comparison, the Markovian-embedding approach has two advantages:
  (i) the construction of the Markovian process $\{\zeta_{t}\}_{t=0}^\infty$
  relies on an ad hoc analysis of the uncertainties in the problem, and
  (ii) the dimension of the resulting \emph{augmented state variable}
  $\widetilde{x}_t := (x_t, \zeta_t) \in \mathbb{R}^{n_x} \times Z$ is fixed.
}

We now can rewrite \eqref{eqn:orig} as
\begin{subequations}\label{eqn:mark}
  \begin{align}
    \label{eqn:mark-obj}\min_{\substack{\{u_t(\cdot)\}_{t=0}^\infty}}\;
    & \mathbb{E}_{\bzeta}\left[\sum_{t=0}^\infty \gamma^t \ell(z_t(\bzeta_{0:t}); \zeta_t)\right]\\
    \st\;
    \label{eqn:mark-con-init}&x_{0}(\zeta_0)\sim \mathcal{X}\\
    \label{eqn:mark-con-dyn}&x_{t}=f(z_{t-1}(\bzeta_{0:t-1}); \zeta_t),\; a.s.,\;t\in\mathbb{I}_{\geq 1}\\
    \label{eqn:mark-con-ineq}&g(z_{t}(\bzeta_{0:t}); \zeta_t) \geq 0,\; a.s.,\;t\in\mathbb{I}_{\geq 0},
  \end{align}
\end{subequations}
by redefining $\ell(\cdot,\cdot)\leftarrow \ell(\cdot,\psi(\cdot))$,
$f(\cdot,\cdot)\leftarrow f(\cdot,\psi(\cdot))$, and so on, and by
parameterizing $z_t(\cdot)$ by $\bzeta_{0:t}$ instead of
$\bxi_{0:t}$.

The key benefit of replacing the nonstationary process $\{\xi_t\}_{t=0}^\infty$ with the higher-dimensional Markov process $\{\zeta_t\}_{t=0}^\infty$ is that we can express the problem in the fixed policy space \FP{parameterized by the augmented state $\widetilde{x}_t = (x_t, \zeta_t)$.
  The augmented state $\widetilde{x}_t$ with control $u_t$ defines a Markov decision process, and it is well known that infinite-horizon Markov decision processes admit a fixed policy solution~\cite{bertsekas2012dynamic}.
Thus, no optimality gap is introduced by the reformulation of
Problem~\eqref{eqn:mark} into}
\begin{subequations}\label{eqn:pol}
  \begin{align}
    \min_{\pi(\cdot)}\;
    & \eqref{eqn:mark-obj}\\
    \st\;
    &\eqref{eqn:mark-con-init} - \eqref{eqn:mark-con-ineq}\label{eqn:pol-con}\\
    &u_t(\bzeta_{0:t}) = \pi (x_t(\bzeta_{0:t}),\zeta_t),\; t\in\mathbb{I}_{\geq 0}.\label{eqn:pol-pol}
  \end{align}
\end{subequations}
Instead of optimizing w.r.t. the control process $\{u_t(\cdot)\}_{t=0}^{\infty}$,
Problem~\eqref{eqn:pol} looks for a solution as a fixed control policy
$\pi:\mathbb{R}^{n_x}\times \mathbb{R}^{n_\zeta}\rightarrow \mathbb{R}^{n_u}$.

In large-scale applications where, for instance, the high-dimensional
process represents weather or climate variables, $\bzeta_{0:t}$ is
typically high dimensional. Thus, its direct computation may be
impractical. To alleviate this issue, one can use a variable selection procedure
to reduce the dimension of the problem to essential variables that
are evolved directly. In contrast, the effect of the latent variables
is approximated. Moreover, the predictability window is limited if the
dynamic system is chaotic. However, some variables have longer
predictability widows than do others. The predictability can be increased
by coarsening the spatial or temporal resolutions. 

\section{Approximations}
Although the problem in \eqref{eqn:pol} is still infinite-dimensional (decision variable is a function whose domain is continuous, the expectation is taken with respect to continuous random variables, and the time horizon is infinite), \eqref{eqn:orig} can be approximated as a finite-dimensional nonlinear program by applying several approximation strategies.
In particular, we apply  function approximation,  sample average approximation (via deterministic sampling), and  temporal truncation. While these approximation strategies are not fully rigorous, we will justify why  they can yield a near-optimal solution.

\subsection{Finite-Dimensional Nonlinear Programming Formulation}
We are now set to state the approximate problem:
\begin{subequations}\label{eqn:approx}
  \begin{align}
    \min_{\theta \in \mathbb{R}^{P}}\;
    & \sum_{i=1}^S\sum_{t=0}^T \gamma^t \ell(z^{(s)}_t(\bzeta^{(s)}_{0:t}); \zeta^{(s)}_t)\\
    \st\;
    &x^{(s)}_{0} = \overline{x}^{(s)},\;  s\in\mathbb{I}_{[1,S]}\\
    &x^{(s)}_{t}=f(z^{(s)}_{t-1}(\bzeta^{(s)}_{0:t}); \zeta^{(s)}_t),\; t\in\mathbb{I}_{[0,T]}, s\in\mathbb{I}_{[1,S]}\\
    &g(z^{(s)}_{t}(\bzeta^{(s)}_{0:t}); \zeta^{(s)}_t) \geq 0,\; t\in\mathbb{I}_{[1,T]}, s\in\mathbb{I}_{[1,S]}\\
    &u^{(s)}_t(\bzeta^{(s)}_{0:t}) = \pi_\theta(x^{(s)}_t,\zeta^{(s)}_t),\; t\in\mathbb{I}_{[1,T]}, s\in\mathbb{I}_{[1,S]} .
  \end{align}
\end{subequations}
The policy $\pi(\cdot)$ is assumed to be approximated by a certain function approximator $\pi_\theta(\cdot)$ with parameter $\theta$.
The samples are obtained by drawing $S$ samples $\{(\overline{x}^{(s)},\bzeta^{(s)}_{0:T})\}_{s=1}^S$ from the known distribution of the uncertainty.
Using these samples, one can replace the expectation in \eqref{eqn:orig-obj}  by the sample average, and the almost-sure condition in the constraints \eqref{eqn:orig-dyn}-\eqref{eqn:orig-con} can be replaced by the constraints enforced over the finite sample points.
The temporal truncation is performed by replacing the infinite-horizon summation in \eqref{eqn:pol} with a finite-horizon summation over $\mathbb{I}_{[0,T]}$ and enforcing dynamic and algebraic constraints only over $\mathbb{I}_{[0,T]}$.
Now the only decision variables in \eqref{eqn:approx} are the finite-dimensional policy parameters $\theta$.
Further, the number of constraints is also finite, and the objective and the constraint functions do not involve infinite-dimensional operators.
Still, the problem is computationally challenging to solve because of the potentially large number of objective terms and constraints.

Applying approximation strategies to \eqref{eqn:pol} rather than \eqref{eqn:orig} has several benefits.
First, one does not need to solve the problem online; if the temporal truncation strategy is applied to \eqref{eqn:orig}, the problem needs to be resolved every several time intervals (called receding-horizon control), since the solution of the approximated problem does not cover the entire infinite temporal horizon.
Second, the target function to be approximated becomes much simpler in \eqref{eqn:pol} than \eqref{eqn:orig} thanks to the time-independent nature of the fixed policy $\pi(\cdot)$.

At the same time, however, properly approximating \eqref{eqn:pol} is likely to be also computationally expensive.
For instance, it is well known that directly constructing the optimal policy based on active set analysis is intractable even for medium-size linear-quadratic problems because of the combinatorial explosion of possible active sets \cite{bemporad2002explicit}.
Thus, the policy representation architectures, such as neural networks or nonlinear basis functions, need to have good approximating properties for the resulting control of \eqref{eqn:mark}.

\subsection{Discussions}
A natural question in the context of \eqref{eqn:approx} is: {\it Does the minimum of \eqref{eqn:approx} approach to the minimum of \eqref{eqn:pol} as $S,T,P\rightarrow \infty$?}
Addressing this question is nontrivial, and establishing a rigorous end-to-end near-optimality guarantee is extremely challenging.
Here we provide some justification based on the existing theoretical results as well as indicate some open problems for the resulting method.

First, it is well known that certain function approximation paradigms, such as neural networks, enjoy the universal approximation property \cite{cybenko1989approximation}.
Thus, by using neural networks with an adequate number of hidden nodes, one may expect that the approximation quality improves.
In this context, it is of interest to characterize how many hidden nodes one needs in order to make the approximation accurate up to the given error bound.

Second, several works in the literature have studied the sample average approximation strategies for stochastic nonlinear programs.
In particular, several works have established asymptotic convergence \cite{shapiro1990differential,shapiro1993asymptotic}, convergence rates \cite{shapiro2000rate}, and non-asymptotic guarantees \cite{oliveira2022sample}.
However, the formulation studied in these works does not cover the almost-sure constraint in \eqref{eqn:pol}.
Thus, the analysis does not directly apply to our setting, and accordingly it is of interest to establish the asymptotic/non-asymptotic guarantee for the sample average approximation of stochastic NLPs with almost-sure constraint.

The effect of temporal truncation has been studied in several recent works \cite{na2020exponential,shin2022exponential,lin2021perturbation}.
The key idea is that dynamic optimization problems exhibit the exponential decay of sensitivity under regularity conditions given by the stabilizability and detectability, and this leads to the near-optimality of the policies with truncated prediction horizons.
This property is  studied primarily under linear-quadratic settings, however, and the effect of temporal truncation under a nonlinear, constrained setting is still not well understood.
Thus, the effect of temporal truncations in different settings needs to be studied further.

Several alternative approaches exist for approximately solving \eqref{eqn:pol}.
Approximating the sequential decision policies in a learning paradigm has made great progress, culminating in the success of AlphaGo \cite{silver2017mastering}.
However, it cannot handle nonstationary exogenous signals \cite{bertsekas2012dynamic}, and the standard algorithms are not suitable for rigorously handling constraints.
Further, stochastic algorithms for various machine learning problems (e.g., stochastic gradient descent and stochastic quasi-Newton methods \cite{bottou2018optimization}) have been studied for solving stochastic programs appearing in machine learning problems, but they are also underdeveloped for handling general nonlinear constraints.

\section{Algorithms}
We are now set to discuss the numerical solution strategy for \eqref{eqn:approx}.
We first write \eqref{eqn:approx} in a more abstract form:
\begin{align}\label{eqn:comp}
  \min_{\theta \in \mathbb{R}^{P}}\; \sum_{i=1}^S L^{(i)}(\theta) \; \st\;G^{(i)}(\theta) \geq 0,\; i=1,\cdots,S.
\end{align}
Recall that \eqref{eqn:comp} is a minimization problem of the summation of the sampled stage costs over the truncated horizon, subject to the algebraic constraints for each time stage and each sample.
Further note that the temporal structure is collapsed into the expressions $L^{(i)}(\cdot)$ and $G^{(i)}(\cdot)$; thus, the expressions in $L^{(i)}(\cdot)$ and $G^{(i)}(\cdot)$ embed the $T$-step forward simulation with the policy $\theta$ and dynamics \eqref{eqn:dyn}.
Computing the solution of \eqref{eqn:comp} can be performed by exploiting the structure.

\subsection{Automatic Differentiation (AD)}
Applying a Newton or quasi-Newton method for \eqref{eqn:comp} requires evaluating the first- and/or second-order derivatives of \eqref{eqn:comp}.
Using efficient AD is important, since evaluating the derivatives of the expressions in \eqref{eqn:comp} can be  expensive because they can embed long-horizon dynamical simulations.
One can observe that $L^{(i)}(\cdot)$ and $G^{(i)}(\cdot)$ for each $i$ are independent. Thus the associated derivative evaluations can be performed in parallel.
Moreover, for all $i$, both $L^{(i)}(\cdot)$ and $G^{(i)}(\cdot)$ share the same structure in their expressions and differ only  in the numerical values of the sampled data $\zeta^{(i)}$.
This means that the evaluation of first- and second-order derivatives can be performed by using SIMD parallelism, which is suitable for running on modern GPUs.
In particular, by evaluating the objective gradient and constraint Jacobian via batched operations, significant speed-up can be achieved.
The implementation of the batched first- and second-order derivative evaluations using  state-of-the-art AD tools such as Zygote.jl \cite{Zygote.jl-2018}, JAX~\cite{bradbury2018jax}, and PyTorch \cite{paszke2017automatic} will be investigated in the future.

\subsection{Condensed-Space Interior-Point Method }
The problem in \eqref{eqn:comp} can be solved with the standard constrained nonlinear optimization algorithms, such as the augmented Lagrangian method, sequential quadratic programming, and IPM \cite{nocedal1999numerical}.
Here we describe  IPM, but using other algorithms is also possible.

The standard IPM introduces the slack variables and replaces the inequality constraints with log-barrier functions:
\begin{align*}
  \min_{\theta , s}\;& \sum_{i=1}^S L^{(i)}(\theta) - \mu\log s^{(i)} \\
  \st\;&G^{(i)}(\theta) = s^{(i)},\; \forall i=1,\cdots,S.
\end{align*}
Here $\mu>0$ is the penalty parameter, and $s$ is the slack variable.
Then, the primal-dual search direction is computed by applying the regularized Newton steps (approximated Newton steps in the case of quasi-Newton) for the KKT conditions of the barrier subproblem.
This can be performed by solving a linear system of the following structure:

{\scriptsize
  \begin{align}\label{eqn:newton}
  \begin{bmatrix}
    \sum_{i=0}^S H^{(i)} + \Sigma_\theta & &&& J^{{(1)}\top}& \cdots & J^{{(S)}\top}\\
    &\Sigma_1&&&I\\
    &&\ddots&&&\ddots\\
    &&&\Sigma_S&&&I\\
    J^{(1)}& I\\
    \vdots&&\ddots\\
    J^{(S)}&&&I
  \end{bmatrix}
  \begin{bmatrix}
    d_\theta\\
    d^s_1\\
    \vdots\\
    d^s_1\\
    d^\lambda_1\\
    \vdots\\
    d^\lambda_S
  \end{bmatrix}=
  \begin{bmatrix}
    r_\theta\\
    r^s_1\\
    \vdots\\
    r^s_1\\
    r^\lambda_1\\
    \vdots\\
    r^\lambda_S
  \end{bmatrix}.
  \end{align}
  \normalsize
  Here, $H^{(i)}$ is the exact or approximate Hessian of Lagrangian $\nabla_{\theta\theta}^2 L^{(i)}(\theta) + \lambda^{(i)} G^{(i)}(\theta)$ for sample $i$; $J^{(i)}$ is the Jacobian of the constraint $\nabla_{\theta}G^{(i)}(\theta)$ for sample $i$; $\Sigma_{\theta}$ and $\Sigma_{1},\cdots,\Sigma_{S}$ are the terms that derive from regularization and the barrier term; $d_\theta,d^s_1,\cdots,d^s_S,d^\lambda_1,\cdots,d^\lambda_S$ are step directions; and $r_\theta,r^s_1,\cdots,r^s_S,r^\lambda_1,\cdots,r^\lambda_S$ are the right-hand  side terms that derive from the KKT residual.
  The step direction obtained by solving \eqref{eqn:newton} can be used for updating the interior-point primal-dual iterate, with the help of the line search method; and by repeating the procedure with the decreasing sequence of penalty parameters $\mu\searrow 0$, one can converge to the stationary point (which often happens to be the local minimum) of \eqref{eqn:comp}.
  For the details of the IPMs, the readers are pointed to \cite{nocedal1999numerical}.
  }

  Directly solving the KKT system in \eqref{eqn:newton} is computationally intractable since the Hessian and Jacobian matrices are dense and the size of the system scales linearly with $T$ and $S$.
  To avoid factorizing this large matrix, we apply the {\it condensation} technique, based on Schur complement decomposition.
  In particular, by exploiting that $\Sigma_1, \cdots, \Sigma_S\succ 0$ (always holding due to the nature of IPM), one can solve the system in \eqref{eqn:newton} by solving
\begin{align}\label{eqn:condense}
  \Bigg(\sum_{i=0}^S H^{(i)}& + \Sigma_\theta + \sum_{i=0}^S J^{(i)\top} \Sigma_s J^{(i)}\Bigg) d_\theta  \\
  &= r_\theta  - \sum_{i=0}^S J^{(i)\top} (r^s_i - \Sigma_i r^\lambda_i),\nonumber
\end{align}
and then computing the rest of the directions based on \eqref{eqn:newton}.
Solving the system in \eqref{eqn:condense} is significantly cheaper than directly solving \eqref{eqn:newton}, since the system size is much reduced, and the system now becomes a positive definite system (after regularization), as opposed to an indefinite system.
This allows using the Cholesky factorization, which is highly scalable on GPU/SIMD architectures.

For  efficient implementation of the condensation strategy in \eqref{eqn:condense}, one needs an implementation of an implicit Jacobian operator that does not explicitly evaluate the full dense Jacobian but has the capability to perform Jacobian vector product and compute $J^{(i)\top}\Sigma^{(i)} J^{(i)}$.
This allows us to significantly reduce the memory footprint and removes the potential memory bottleneck.
Further, applying this method requires an efficient implementation of the IPM that ports the solution of the dense KKT systems on the GPU.
Currently, MadNLP.jl \cite{shin2020graph} has the capability to run majority of the computation on GPU, and several additional features (the use of approximate Hessian via BFGS algorithm \cite{nocedal1999numerical}) are under development.

\section{Numerical Experiments}
We demonstrate the proposed approach with the following example, inspired by the heating/cooling problem under time-varying and uncertain energy prices and disturbances:
\begin{align*}
  \ell(x,u;\xi)&=10^{-3}\|x\|^2 + \|u\|^2 +
                 \xi\begin{bmatrix}
                   0.3&0.3&0.3
                 \end{bmatrix}^\top
                 u \\
  g(x,u;\xi)&=\begin{bmatrix}
    0.9&-0.05&0\\
    -0.05&0.9&-0.05\\
    0&-0.05&0.9&\\
  \end{bmatrix}x + u +
  \begin{bmatrix}
    0.1\\
    0.1\\
    0.1
  \end{bmatrix}
  \xi\\
  g(x,u;\xi)&\geq 0\iff\\
  \begin{bmatrix}
    -0.2\\
    -0.2\\
    -0.2\\
  \end{bmatrix}&\leq x\leq  \begin{bmatrix}
    0.2\\
    0.2\\
    0.2\\
  \end{bmatrix},\;
  \begin{bmatrix}
    -0.03\\
    -0.03\\
    -0.03\\
  \end{bmatrix}
  \leq u\leq  \begin{bmatrix}
    0.03\\
    0.03\\
    0.03
  \end{bmatrix},
\end{align*}
where $\{\xi_t\}_{t=0}^\infty$ satisfies
\begin{align*}
    \zeta_{t}
    &=
      \begin{bmatrix}
        0.955&  0.295&0\\
        -0.295&  0.955&0\\
        0& 0 & 0 \\
      \end{bmatrix}\zeta_{t-1}  +
  \begin{bmatrix}
    0\\0\\1
  \end{bmatrix}w_t
      ,\; t=1,2,\cdots,\\
    \xi_t
    &=
    \begin{bmatrix}
      1 & 0& 1\\
    \end{bmatrix}\zeta_t, \;t=0,1,\cdots,
\end{align*}
for an iid random variable $w_t$.

\subsection{Methods}
We compare three control policies:  the policies obtained by policy optimization (PO),  LQR, and  MPC with a horizon length $20$.
The gain matrix for LQR is computed by solving the algebraic Riccati equation on the augmented state-space model.
We assume that the actuator is saturated when the input constraints are violated.
The MPC problem is formulated by using the standard formulation over the augmented state-space model.
The MPC policy falls back to the LQR policy when the problem becomes infeasible.

The PO problem is formulated with a neural network policy approximator with two hidden layers, with $6$ nodes in each layer, and with $\tanh$ activation functions.
The policies are obtained for a varied number of samples $S$ and the time horizon length $T$.
The validation is performed for 10 sample trajectories with $100$ time stages.
The training data and the validation date are not the same but are drawn from the same distribution.
We performed two sets of experiment.
In the first set, we set $w_t$ to be constantly zero (thus, the system is noise-free);  in the second set, we assume that $w_t$ is drawn from the uniform distribution over $[-1,1]$.
The code is implemented in the Julia language with the help of several packages: ForwardDiff.jl (automatic differentiation) \cite{RevelsLubinPapamarkou2016}, NLPModels.jl (nonlinear programming abstraction) \cite{orban-siqueira-nlpmodels-2020}, and MadNLP.jl (solving policy optimization problems) \cite{shin2020graph}.
The results can be reproduced with the code available at \cite{code}.
The policy optimization problems are solved in a reasonable time (all less than $30$ seconds) on a laptop computer with an Intel i9-11900H CPU.

\subsection{Results}
\begin{figure*}
  \centering
  \includegraphics[width=.48\textwidth]{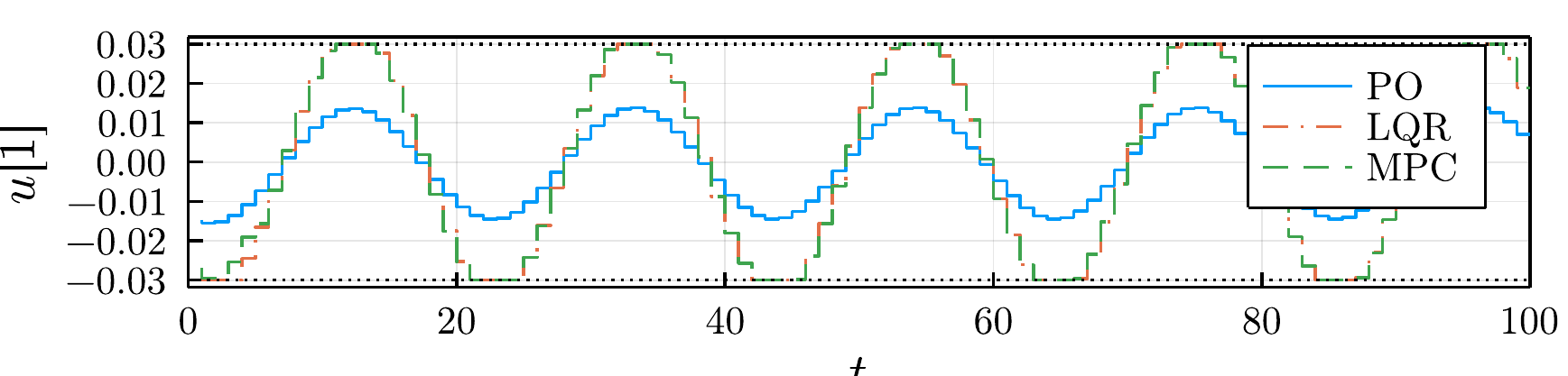}
  \includegraphics[width=.48\textwidth]{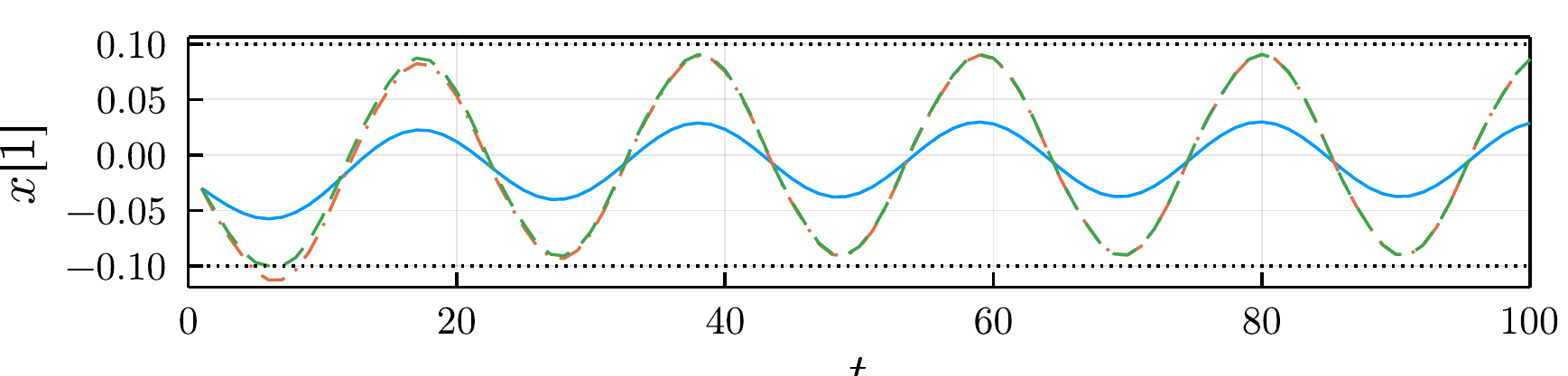}\\
  \includegraphics[width=.48\textwidth]{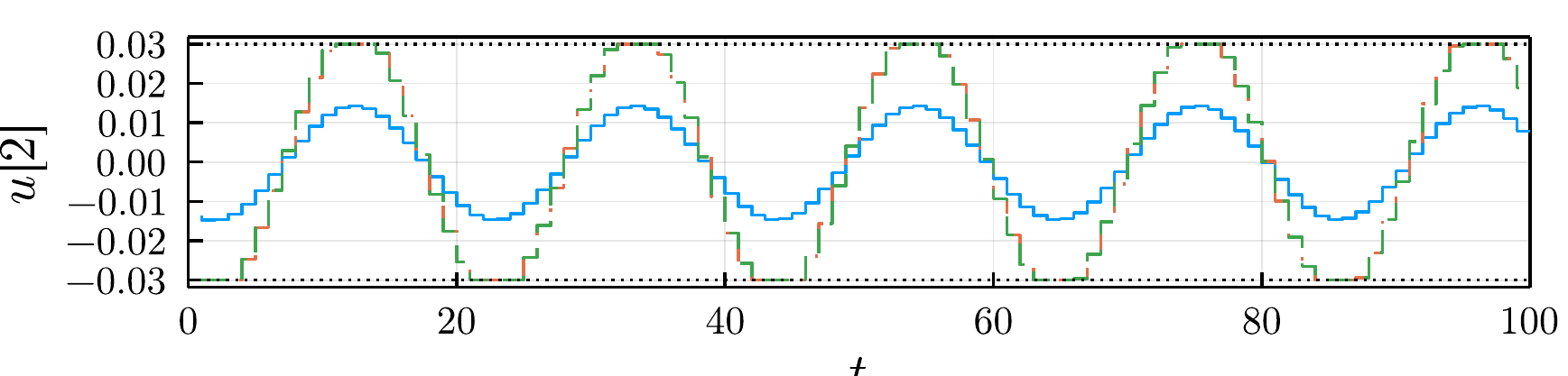}
  \includegraphics[width=.48\textwidth]{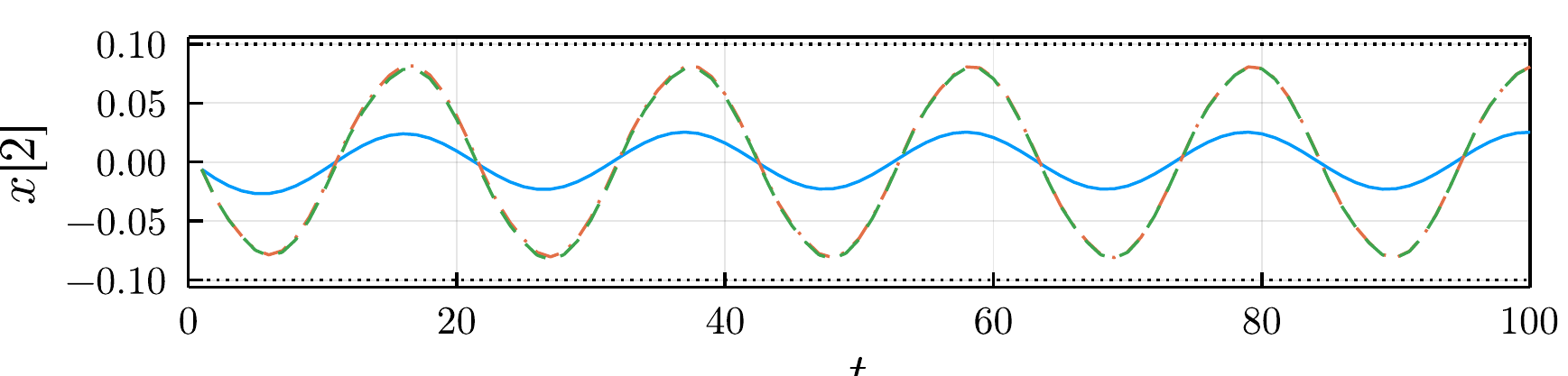}\\
  \includegraphics[width=.48\textwidth]{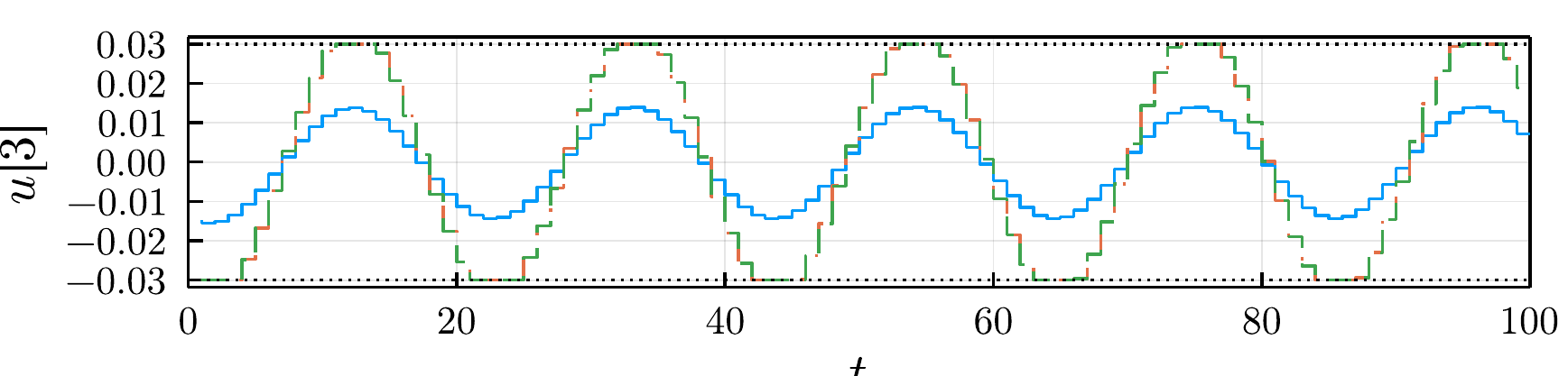}
  \includegraphics[width=.48\textwidth]{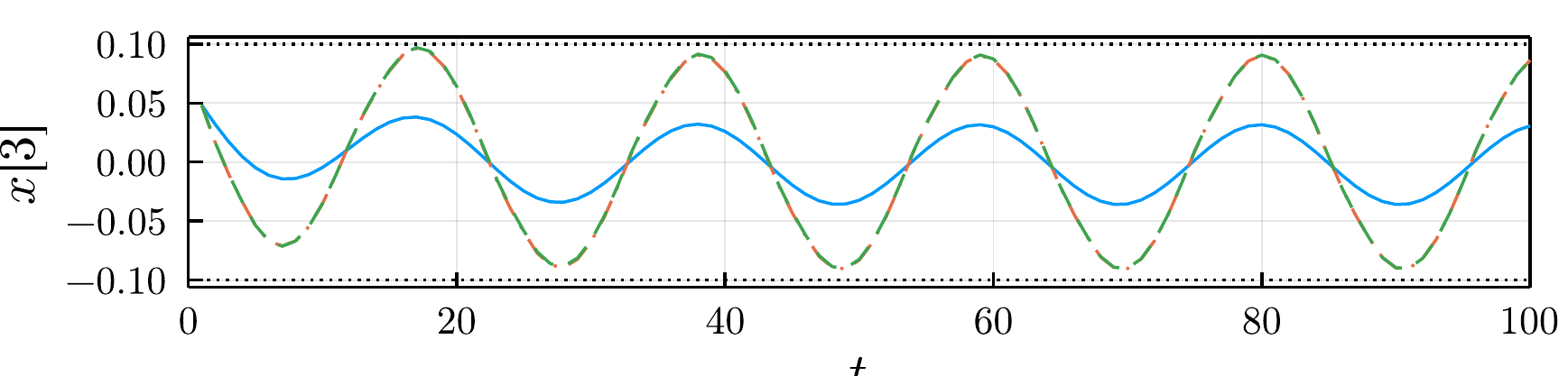}
  \caption{Closed-loop simulation of PO, LQR, and MPC (nominal).}\label{fig:nominal}
  \vspace{.05in}
  \includegraphics[width=.48\textwidth]{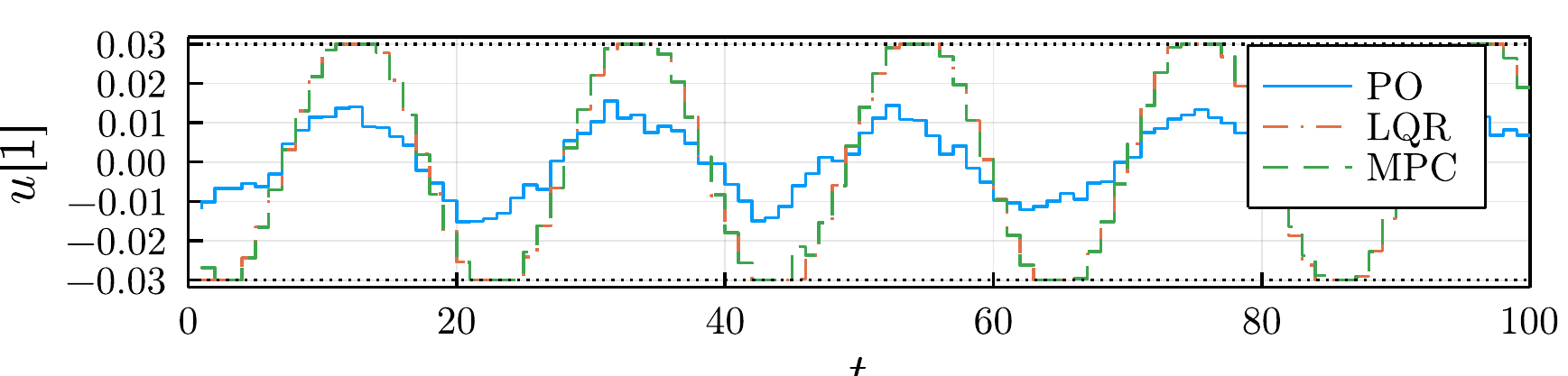}
  \includegraphics[width=.48\textwidth]{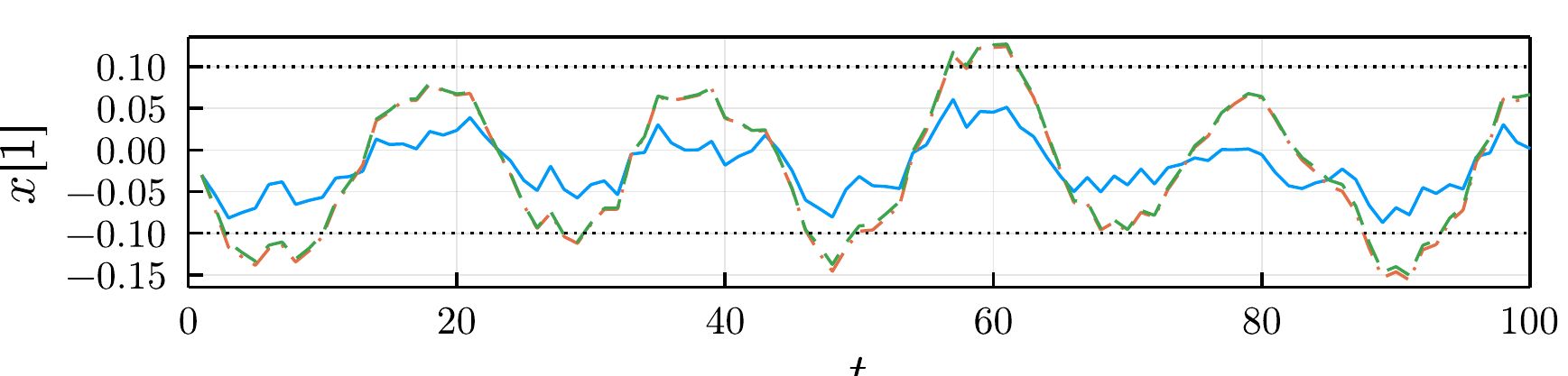}\\
  \includegraphics[width=.48\textwidth]{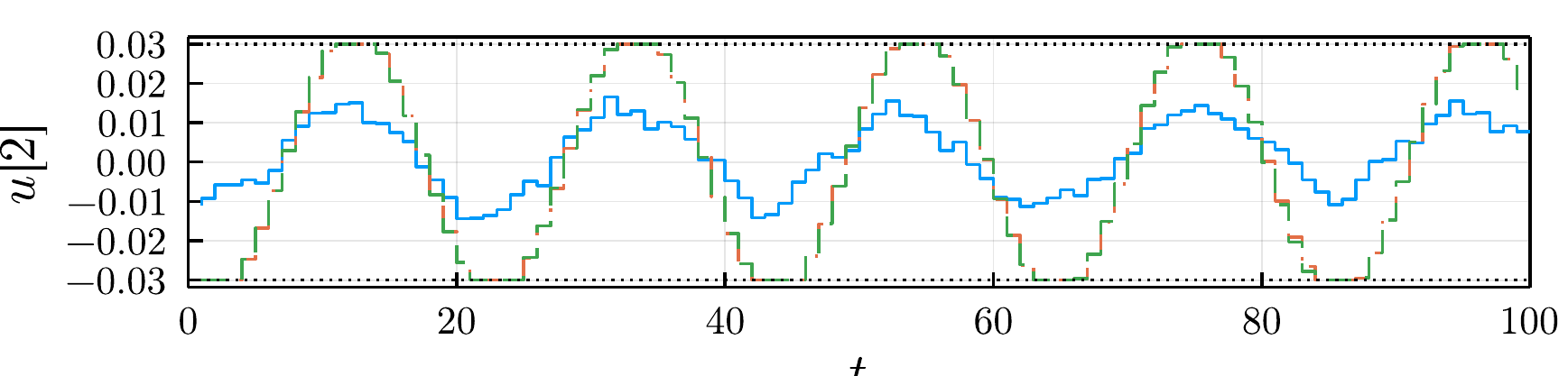}
  \includegraphics[width=.48\textwidth]{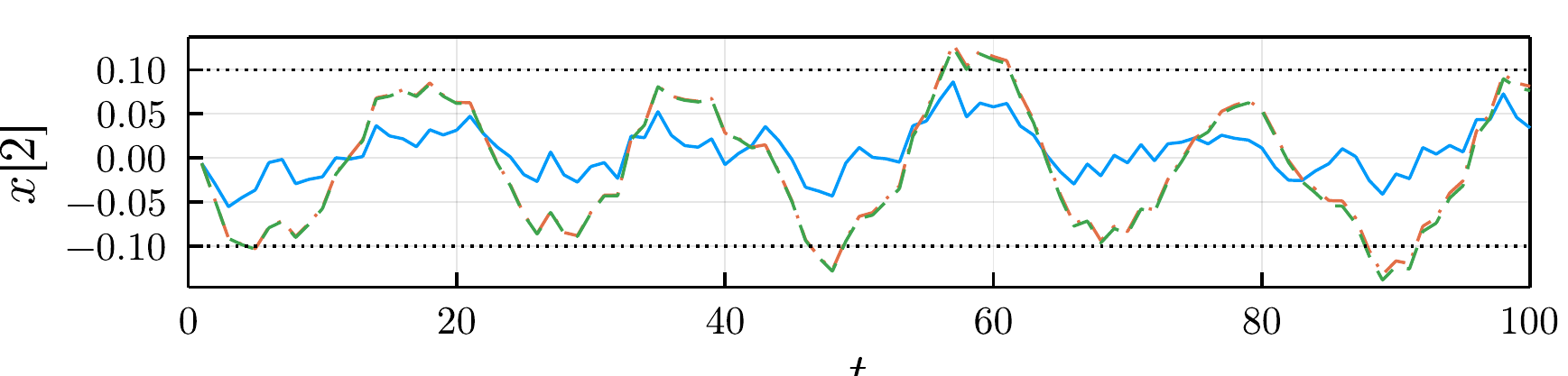}\\
  \includegraphics[width=.48\textwidth]{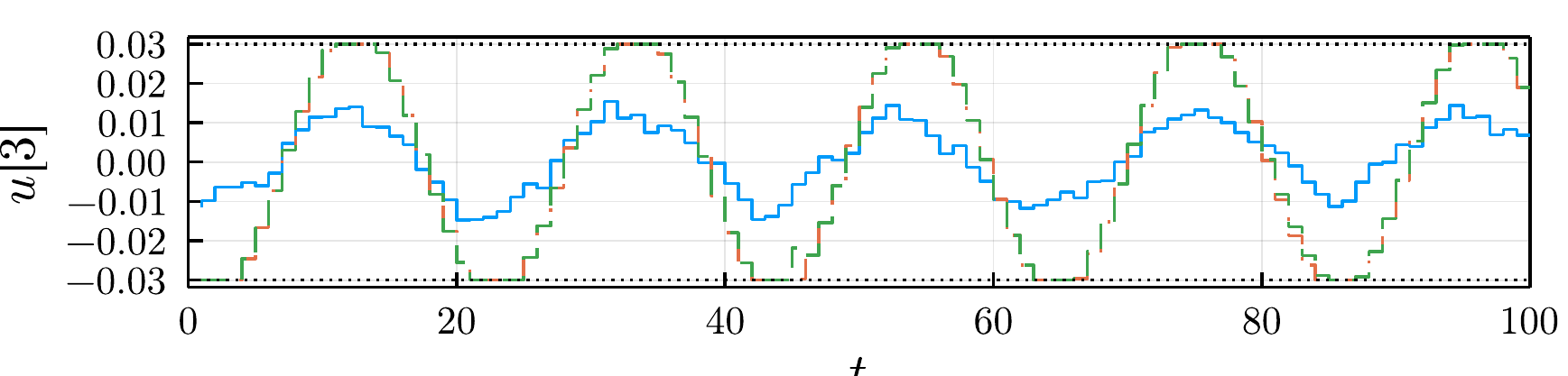}
  \includegraphics[width=.48\textwidth]{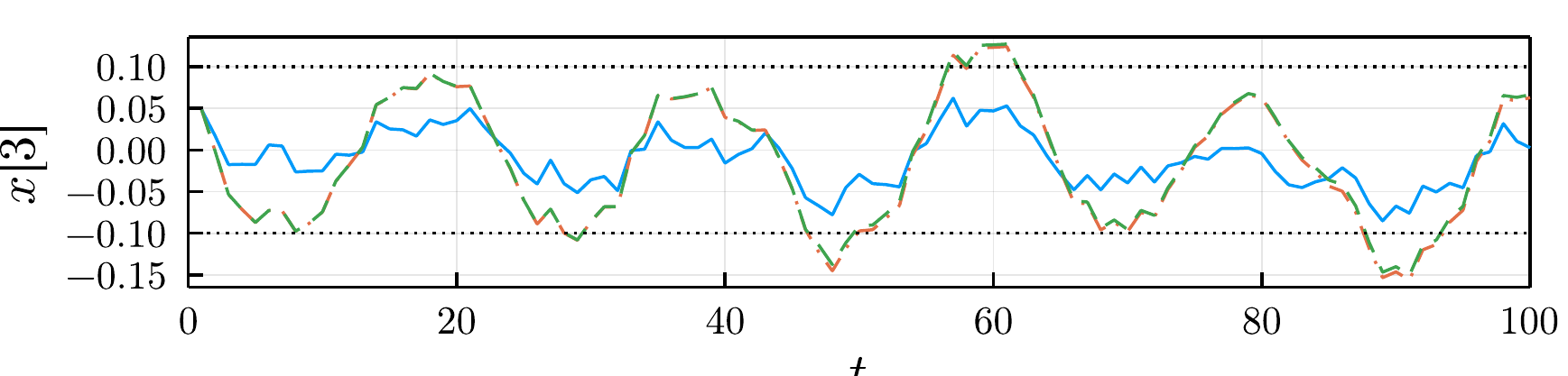}
  \caption{Closed-loop simulation of PO, LQR, and MPC (noisy).}\label{fig:noisy}
\end{figure*}

First we discuss the noise-free case.
The performance and constraint violation results are shown in Table \ref{tbl:nominal}, and one of the closed-loop controller profiles for a representative scenario is shown in Figure \ref{fig:nominal}.
In Table \ref{tbl:nominal}, performance refers to the sample-averaged closed-loop performance index over the validation data, and the constraint violation represents the average number of time points where the constraints are violated for every 100 time stages on the validation data.
One can see that LQR and MPC perform approximately the same, but LQR violates the constraints for several time points.
Further, one can see that the policies obtained from PO with high $S$ and $T$ do not violate the constraints.
However, the closed-loop performance indexes were significantly worse than for LQR and MPC.

\begin{table*}[t]
  \centering
  \caption{Performance Comparison of PO, LQR, and MPC (nominal)}
  \begin{tabular}{|c|c|c|c|c|c|c|c|c|c|c|c|c|c|}
    \hline
    &\multicolumn{5}{c|}{PO} & \multirow{2}{*}{LQR} & \multirow{2}{*}{MPC} \\
    \cline{2-6}
    &
                 $S\backslash T$&$5$&$10$&$15$&$20$&&\\
    \hline
    \multirow{4}{*}{\shortstack{Performance\\(the lower the better)}}
    &$5$& 0.000535093& 0.00206096& -0.0660549& -0.0660536& \multirow{4}{*}{-0.07679554801271361}& \multirow{4}{*}{-0.0766397836932208}\\
    \cline{2-6}
    &$10$& -0.0657874& 0.000662422& -0.0561135& -0.0536368&&\\
    \cline{2-6}
    &$15$& 0.00374853& 0.00131267& -0.0544788& -0.0531895&&\\
    \cline{2-6}
    &$20$& 0.00118073& 0.000250603& -0.0533167& -0.052409&&\\
    \hline
    \multirow{4}{*}{\shortstack{Constraint Violation\\(the lower the better)}}
    &$5$ & 0.0& 32.5& 6.66& 6.48& \multirow{4}{*}{5.9} & \multirow{4}{*}{0.0}
    \\
    \cline{2-6}
    &$10$ & 7.59& 0.7& 0.31& 0.04&&\\
    \cline{2-6}
    &$15$ & 79.46& 0.23& 0.0& 0.0&&\\
    \cline{2-6}
    &$20$ & 0.0& 0.0& 0.0& 0.0&&\\
    \hline
  \end{tabular}
  \vspace{.05in}
  \label{tbl:nominal}
  \caption{Performance Comparison of PO, LQR, and MPC (noisy)}
  \begin{tabular}{|c|c|c|c|c|c|c|c|c|c|c|c|c|c|}
    \hline
    &\multicolumn{5}{c|}{PO} & \multirow{2}{*}{LQR} & \multirow{2}{*}{MPC} \\
    \cline{2-6}
    &
                 $S\backslash T$&$5$&$10$&$15$&$20$&&\\
    \hline
    \multirow{4}{*}{\shortstack{Performance\\(the lower the better)}}
    &$5$& 0.000184756& 0.000351228& 0.000713321& -0.0691453& \multirow{4}{*}{-0.07755931501730873}& \multirow{4}{*}{ -0.07734653866399706}\\
    \cline{2-6}
    &$10$& -0.0661179& 0.000192822& -0.0540974& -0.0573766&&\\
    \cline{2-6}
    &$15$& 0.000607659& 0.00288476& -0.0536313& -0.0558113&&\\
    \cline{2-6}
    &$20$& 0.000461463& -0.0480934& -0.054308& -0.0476584&&\\
    \hline
    \multirow{4}{*}{\shortstack{Constraint Violation\\(the lower the better)}}
    &$5$ & 1.24& 1.26& 2.86& 27.8& \multirow{4}{*}{16.75} & \multirow{4}{*}{14.39}
    \\
    \cline{2-6}
    &$10$ & 10.77& 1.93& 16.59& 4.09&&\\
    \cline{2-6}
    &$15$ & 1.35& 12.34& 2.74& 1.36&&\\
    \cline{2-6}
    &$20$ & 3.47& 8.51& 1.05& 0.45&&\\
    \hline
  \end{tabular}
  \label{tbl:noisy}
\end{table*}

We now discuss the noisy case.
The performance and constraint violation results are shown in Table \ref{tbl:noisy}, and one of the closed-loop controller profiles for a representative scenario is shown in Figure \ref{fig:noisy}. While MPC has a recursive feasibility guarantee for nominal settings \cite{rawlings2017model}, it may violate the constraints if the system is subject to uncertainty.
In particular, the event sequence and dynamics in Section \ref{sec:setting-1} prohibit MPC from accurately predicting the next state.
The results in Table \ref{tbl:noisy} and Figure \ref{fig:noisy} reveal this limitation of MPC; the constraints are violated for many time stages.
Also, LQR violates the constraints for many time stages because of its inherent limitation in handling the constraints as well as the effect of uncertainties.
However, the policy obtained by PO with high $T$ and $S$ violates the constraints much less frequently.
On the other hand, the closed-loop performance index of the policies obtained by PO is still much worse than that of MPC and LQR.

These results suggest that a control policy with respectable performance can be obtained through the PO approach with a moderate number of samples and time horizon length.
However, the closed-loop performance index was found to be significantly worse than that of the LQR and MPC.
For the noisy setting, the optimal solution may indeed perform significantly worse than MPC and LQR---since they always satisfy the constraints, they make conservative moves---but the results for the noise-free setting suggest that the policy obtained by PO may be more conservative than what is necessary.
This limitation might have been caused by the (i) inadequate number of parameters in the policy approximator, (ii) the inadequate number of samples, and (iii) the short time horizon.
In the future, the capability of the PO approaches will be investigated further with more scalable implementation and realistic case studies.
The method also will be compared with different methods in stochastic model predictive control \cite{mesbah2016stochastic,lucia2020stability,bernardini2009scenario}, which are expected to perform better than the deterministic MPC considered here.

\section{Conclusions}

We have presented a constrained policy optimization approach for solving stochastic control problems under nonstationary uncertainty based on state augmentation, reformulation into policy optimization, and finite-dimensional approximation.
The approximation involves function approximation, sample average approximation, and temporal truncation, each of which can be justified to some extent, but the precise error analysis remains  an open question.
Even after the approximation, the solution of the approximate problem is nontrivial because of the complex function expressions and the large number of samples.
We have highlighted several computational challenges and proposed scalable computing strategies that exploit the separable structure and are suitable for running on SIMD/GPU architecture.
The numerical results suggest that the policies obtained by the policy optimization method can achieve a respectable control performance while having a better capability of satisfying the constraints compared with LQR and MPC.

\bibliographystyle{IEEEtran}
\bibliography{ref}

\vspace{0.1cm}
\begin{flushright}
	\scriptsize \framebox{\parbox{2.5in}{Government License: The
			submitted manuscript has been created by UChicago Argonne,
			LLC, Operator of Argonne National Laboratory (``Argonne").
			Argonne, a U.S. Department of Energy Office of Science
			laboratory, is operated under Contract
			No. DE-AC02-06CH11357.  The U.S. Government retains for
			itself, and others acting on its behalf, a paid-up
			nonexclusive, irrevocable worldwide license in said
			article to reproduce, prepare derivative works, distribute
			copies to the public, and perform publicly and display
			publicly, by or on behalf of the Government. The Department of Energy will provide public access to these results of federally sponsored research in accordance with the DOE Public Access Plan. http://energy.gov/downloads/doe-public-access-plan. }}
	\normalsize
\end{flushright}

\end{document}

%% file: preamble.tex
\usepackage{amsmath,amssymb,amsfonts} 
\usepackage[linktocpage=true,colorlinks=true,linkcolor=blue,citecolor=blue,urlcolor=blue]{hyperref}
\usepackage{tikz}
\usepackage{color}
\usepackage{ifthen}

\usetikzlibrary{shapes}
\usetikzlibrary{calc}

\newtheorem{assumption}{Assumption}

\providecommand{\keywords}[1]
{
  \small
  \textbf{\textit{Keywords---}} #1
}

\newcommand{\bxi}{\boldsymbol{\xi}}

\newcommand{\bzeta}{\boldsymbol{\zeta}}

\newcommand{\st}{\mathop{\text{\normalfont s.t.}}}

\newcommand{\bbR}{\mathbb{R}}
\newcommand{\bbP}{\mathbb{P}}

\newcommand{\cF}{\mathcal{F}}

\newcommand{\FP}[1]{{\color{black} #1}}

\usepackage{algorithm}
\usepackage{algpseudocode}

\algdef{SE}[SUBALG]{Indent}{EndIndent}{}{\algorithmicend\ }%
\algtext*{Indent}
\algtext*{EndIndent}

\usepackage[utf8]{inputenc} 
\usepackage[T1]{fontenc}    
\usepackage{hyperref}       
\usepackage{url}            
\usepackage{booktabs}       
\usepackage{amsfonts}       
\usepackage{nicefrac}       
\usepackage{microtype}      